\documentclass[final,3p]{elsarticle} 
\usepackage[utf8]{inputenc}

\journal{Journal of Computational Physics}

\usepackage{amssymb}

\usepackage{graphicx}
\usepackage{xcolor}
\usepackage{mathtools}
\usepackage{easyReview} 

\usepackage{hyperref}

\newcommand{\PP}[1][\alpha ,\beta ]{P^{(#1)}}
\newcommand{\hh}[1][\alpha ,\beta ]{h^{(#1)}}
\newcommand{\xx}[1][\alpha ,\beta ]{\xi ^{(#1)}}
\newcommand{\mat}[1]{\mathbf{#1}}

\newcommand{\tM}{\tilde{\mat{M}}}
\newcommand{\Tr}{\operatorname{Tr}}
\AtBeginDocument{\renewcommand{\Im}{\operatorname{Im}}}

\begin{document}
	
\begin{frontmatter}
	
	
	
	\title{The kernel polynomial method based on Jacobi polynomials}
	
	
	\author[IoffeInst]{I.\,O.~Raikov}
	\ead{IORaikov@gmail.com}
	\author[IoffeInst]{Y.\,M.~Beltukov} 
	\ead{ybeltukov@gmail.com}
	
	\affiliation[IoffeInst]{organization={Ioffe Institute},
		addressline={Politechnicheskaya st. 26}, 
		city={St. Petersburg},
		postcode={194021},
		country={Russia}}
	
	\begin{abstract}
	    The kernel polynomial method based on Jacobi polynomials $\PP_n(x)$ is proposed. The optimal-resolution positivity-preserving kernels and the corresponding damping factors are obtained. The results provide a generalization of the Jackson damping factors for arbitrary Jacobi polynomials. For $\alpha =\pm 1/2$, $\beta =\pm 1/2$ (Chebyshev polynomials of the first to fourth kinds), explicit trigonometric expressions for the damping factors are obtained. The resulting algorithm can be easily introduced into existing implementations of the kernel polynomial method.
	\end{abstract}
	
	
	
	\begin{keyword}
		Kernel polynomial method \sep Jacobi polynomials \sep damping factors \sep non-negative kernels
		
		
		
	\end{keyword}
	
\end{frontmatter}

\section{Introduction}
In numerous disciplines, it is often the case that the analysis of a physical problem can be reduced to studying the properties of a linearised system. These scenarios can range from determining the eigenvalues of quantum systems to identifying the eigenmodes of small atomic vibrations. Traditionally, such problems are tackled by diagonalizing the corresponding matrices to extract their eigenvalues. However, given the numerous degrees of freedom in these systems, the diagonalization process typically demands significant computational resources. An alternative approach is the Kernel Polynomial Method (KPM)~\cite{Silver-1994-densities-states-megadimensional,Silver-1996-kernel-polynomial-approximations,Weisse-2006-kernel-polynomial-method}. This method allows the computation of properties such as the density of states and the Green function by decomposing them into orthogonal polynomials, eliminating the need to search for eigenvalues and eigenvectors individually. KPM offers the advantage of analyzing systems with millions of degrees of freedom without the use of supercomputers. It has found numerous applications in contemporary research, including investigating atomic vibrations~\cite{Beltukov-2016-boson-peak-iofferegel}, quantum transport phenomena~\cite{Bertini-2021-finitetemperature-transport-onedimensional}, 
topological insulators~\cite{Varjas-2020-computation-topological-phase}, 
Anderson localization~\cite{Khan-2021-probing-bandcenter-anomalya}, 
non-Hermitian systems~\cite{Chen-2023-topological-spin-excitations}, 
and twistronics~\cite{Carr-2017-twistronics-manipulating-electronic}.
Furthermore, KPM can be effectively utilized in calculating the spectral density of graph Laplacians~\cite{Braverman-2022-sublinear-time-spectral}.

The KPM was described in detail in \cite{Weisse-2006-kernel-polynomial-method}. It is closely related to several other methods, such as the Delta-Gauss-Legendre (DGL) method~\cite{Lin-2016-approximating-spectral-densities} and the Lanczos spectroscopic procedure~\cite{Lin-2016-approximating-spectral-densities}. It has some connection to the Lanczos recursion method -- for example, it is possible to compute KPM moments from the Lanczos outputs~\cite{Chen-2023-spectrum-adaptive-kernel}. Compared to Lanczos, it does not suffer from the loss of orthogonality, and thus typically has better convergence when estimating the distribution of non-extremal eigenvalues or calculating correlation functions~\cite{Weisse-2006-kernel-polynomial-method}.

Usually~\cite{Weisse-2006-kernel-polynomial-method,Lin-2016-approximating-spectral-densities} the KPM is based on the Chebyshev polynomials of the first kind $T_n = \cos(n\varphi )$, where $\varphi =\arccos x$. Such trigonometric form makes it possible to apply the well-known properties of the Fourier transform. However, the KPM can be based on other orthogonal polynomials with a three-term recurrent formula, such as Legendre~\cite{Lin-2016-approximating-spectral-densities,Qin-2010-linearscaling-momentbased-electronic}, Hermite~\cite{Qin-2010-linearscaling-momentbased-electronic,Qin-2011-calculation-binding-energy} or more generally Jacobi~\cite{Qin-2011-calculation-binding-energy} polynomials, which can give a better approximation depending on the properties of the system under consideration~\cite{Wolf-2015-spectral-functions-time}.

For some tasks, it's important that the approximated DOS doesn't experience Gibbs oscillations (potentially causing it to become negative), which requires the use of proper damping factors depending on the polynomials~\cite{Silver-1994-densities-states-megadimensional,Weisse-2006-kernel-polynomial-method,Lin-2016-approximating-spectral-densities}. The choice of damping factors also affects the spectral resolution. The aforementioned generalizations of the KPM, however, either do not make use of any damping factors or use empirically chosen ones.
Methods related to KPM, such as DGL, typically also don't have non-negativity by default~\cite{Lin-2016-approximating-spectral-densities}. By contrast, the Lanczos recursion method and related methods preserve positivity of the DOS~\cite{Lin-2016-approximating-spectral-densities}.

There are several existing software libraries implementing KPM, including those referenced in~\cite{Groth-2014-kwant-software-packagea, Moldovan-2020-pybinding-v0-python}. However, the existing implementations are limited to utilizing Chebyshev polynomials of the first kind for KPM calculations.

In this paper we propose a generalization of the KPM using the Jacobi polynomials. In Section \ref{sec:kpm} we start by describing the classical Kernel Polynomial Method and our generalization of it to Jacobi polynomials. In Section \ref{sec:optimal} we obtain the optimal (best-resolution) non-negative kernel and the explicit form of the corresponding damping factors. In Section \ref{sec:res} we analyze the resolution of the optimal kernel. In Section \ref{sec:special} the special cases corresponding to Chebyshev polynomials are provided. In Section \ref{sec:asymp_opt_kernel} the asymptotics of the optimal kernel for a large degree of approximation are obtained. Finally, in Section \ref{sec:examples} we demonstrate the application of our method to several problems.

\section{The Kernel Polynomial Method}
\label{sec:kpm}

\subsection{Spectral properties}
The main goal of the KPM is the analysis of the spectral properties of some Hermitian matrix $\mat{M}$ of the size $\mathcal{N}\times \mathcal{N}$. This matrix has eigenvalues $\varepsilon _j$ and the corresponding eigenvectors $|v_j\rangle $ such that 
\begin{equation}
    \mat{M}|v_j\rangle  = \varepsilon _j |v_j\rangle .
\end{equation}
The density of states (DOS) of the matrix $\mat{M}$ is
\begin{equation}
    \rho (\varepsilon ) = \frac{1}{\mathcal{N}} \sum_{j=1}^\mathcal{N} \delta (\varepsilon  - \varepsilon _j).  \label{eq:rho}
\end{equation}
The KPM can be used to make a fast polynomial approximation of $\rho (\varepsilon )$ of the form~\cite{Weisse-2006-kernel-polynomial-method}
\begin{equation}\label{eq:rho_kpm_generic}
	\rho _{\textsc{kpm}}(\varepsilon ) = w(\varepsilon ) \sum_{n=0}^{N-1} s_n P_n(\varepsilon ),
\end{equation}
where $P_n(\varepsilon )$ are polynomials orthogonal with weight $w(\varepsilon )$ on the interval $\varepsilon _\text{min} \leq  \varepsilon  \leq  \varepsilon _\text{max}$ which contains all the eigenvalues $\varepsilon _i$, and $N$ is the approximation order. The coefficients $s_n$ that encode the information about the eigenvalues of $\mat{M}$ are obtained using a stochastic recurrent procedure, which will be discussed later.

In the general case, the sum in Eq.~\eqref{eq:rho_kpm_generic} is not equal to zero at $\varepsilon =\varepsilon _\text{min}$ and $\varepsilon =\varepsilon _\text{max}$.  Therefore, the behavior of the approximation $\rho _{\textsc{kpm}}(\varepsilon )$ near these edge points is given by the behavior of $w(\varepsilon )$. In many problems, the edge behavior of $\rho (\varepsilon )$ is a power law,
\begin{align}
	\rho (\varepsilon ) &\sim (\varepsilon _\text{max}-\varepsilon )^\alpha ,    \label{eq:power_law_left_edge}\\
	\rho (\varepsilon ) &\sim (\varepsilon -\varepsilon _\text{min})^\beta .    \label{eq:power_law_right_edge}
\end{align}
For example, electrons in a semiconductor with parabolic dispersion in the conduction band have the DOS $\rho (\varepsilon ) \sim (\varepsilon -\varepsilon _\text{min})^{d/2-1}$, where $d$ is the dimensionality of the system and $\varepsilon _\text{min}$ is the bottom of the conduction band~\cite{Ashcroft-1976-solid-state-physics}. For atomic vibrations, $\varepsilon _\text{min} = 0$ and the eigenvalues are squares of eigenfrequencies, $\varepsilon _i = \omega _i^2$.
For low-frequency acoustic modes the Debye model predicts the vibrational density of states $\rho (\omega ) \sim  \omega ^{d-1}$\cite{Ashcroft-1976-solid-state-physics}.

Without loss of generality, we can map the eigenvalues of $\mat{M}$ to the interval $[-1,1]$ using the affine transform
\begin{equation}
	\tM = \frac{2}{\varepsilon _\text{max} - \varepsilon _\text{min}} \mat{M} - \frac{\varepsilon _\text{max} + \varepsilon _\text{min}}{\varepsilon _\text{max} - \varepsilon _\text{min}}\mat{I},  \label{eq:M_tilde}
\end{equation}
where $I$ is the identity matrix. The eigenvalues of $\tM$ are
\begin{equation}
	x_i = 2 \frac{\varepsilon _i-\varepsilon _{\text{min}}}{\varepsilon _{\text{max}} - \varepsilon _{\text{min}}} - 1,\enspace -1 \leq  x_i \leq  1.
\end{equation}
After obtaining an approximation for the distribution of $x_i$, denoted as $\rho (x)$, a reverse transform can be used to obtain an approximation for $\rho (\varepsilon )$ \cite{Weisse-2006-kernel-polynomial-method}. If the values $\varepsilon _{\text{max}}$ and $\varepsilon _{\text{min}}$ are unknown, their estimation can be obtained accurately and in linear time via the Lanczos recursion procedure~\cite{Lanczos-1950-iteration-method-solution,Weisse-2006-kernel-polynomial-method}.

Depending on the orthogonal polynomial used, different polynomial approximation can be obtained. The Jacobi polynomials $\PP_n(x)$ present a wide class of orthogonal polynomials. They are orthogonal with weight
\begin{equation}
	w(x) = (1-x)^\alpha  (1+x)^\beta ,
\end{equation}
which match the power law given by Eqs.~\eqref{eq:power_law_left_edge}--\eqref{eq:power_law_right_edge}. Depending on the properties of the matrix $\mat{M}$, different values of the parameters $\alpha $ and $\beta $ can be used for better polynomial approximation.

In the classical KPM, the Chebyshev polynomials of the first kind $T_n(x)$ are used. These polynomials are essentially equivalent, up to a normalization constant, to the Jacobi polynomials with $\alpha =\beta =-1/2$. In order to enhance the accuracy of approximating the DOS across different systems, we propose a generalization of the algorithm to include Jacobi polynomials $\PP_n(x)$ with adjustable parameters $\alpha $ and $\beta $.

\subsection{Jacobi polynomials}
Let us first give the main properties of the Jacobi polynomials $\PP_n(x)$ that will be used for the KPM. The Jacobi polynomials have the following recurrence relation
\begin{align}
    \PP_0(x) &= 1, \\
    \PP_1(x) &= a_0 x + b_0, \\
    \PP_{n+1}(x) &= (a_n x + b_n)\PP_n(x) - c_n \PP_{n-1}(x),
\end{align}
where
\begin{align}
    a_n &= \frac{(2n+\alpha +\beta +1)(2n+\alpha +\beta +2)}{2(n+1)(n+\alpha +\beta +1)},\\ 
    b_n &= \frac{(2n+\alpha +\beta +1)(\alpha ^2-\beta ^2)}{2(n+1)(n+\alpha +\beta +1)(2n+\alpha +\beta )},\\
    c_n &= \frac{(n+\alpha )(n+\beta )(2n+\alpha +\beta +2)}{(n+1)(n+\alpha +\beta +1)(2n+\alpha +\beta )}.
\end{align}
The orthogonality relation for the Jacobi polynomials is
\begin{equation}
    \int_{-1}^1 \PP_n(x)\PP_m(x) w(x) \, dx = \hh_n \delta _{nm},
\end{equation}
where
\begin{equation}
	\hh_n = \frac{2^{\alpha +\beta +1}\Gamma (n+\alpha +1)\Gamma (n+\beta +1)}{(2n+\alpha +\beta +1)\Gamma (n+\alpha +\beta +1)n!}.
\end{equation}

\subsection{Stochastic calculation of the moments}
The DOS $\rho (x)$ can be expanded in terms of the Jacobi polynomials
\begin{equation}
    \rho (x) = w(x) \sum_{n=0}^\infty \frac{\mu _n}{\hh_n} \PP_n(x),
\end{equation}
where $\mu _n$ is the $n$th Jacobi moment of $\rho (x)$:
\begin{equation}
    \mu _n = \int_{-1}^1 \rho (x) \PP_n(x) \, dx = \frac{1}{\mathcal{N}} \Tr\left[\PP_n(\tM)\right].  \label{eq:Jmom_exact}
\end{equation}
However, the direct calculation of traces in Eq.~\eqref{eq:Jmom_exact} is inefficient from the numerical point of view. It requires a number of costly matrix-matrix multiplications during the calculation of the matrix polynomial $\PP_n(\tM)$. Also, each matrix-matrix multiplication rapidly increases the number of nonzero elements even if the matrix $\tM$ is highly sparse.

The moments $\mu _n$ can be calculated stochastically as (see \ref{app:avg})
\begin{equation}
    \mu _n = \overline{\left\langle r\middle|\PP_n(\tM)\middle|r\right\rangle },  \label{eq:Jmom_rnd}
\end{equation}
where overline denotes the averaging over the random $\mathcal{N}$-dimensional unit vector $|r\rangle $. For efficient numerical calculation of the moments, $\mathcal{N}$-dimensional vectors
\begin{align}
    |u_n\rangle  = \PP_n(\tM) |r\rangle 
\end{align}
can be calculated using the following recurrence relations:
\begin{align}
    |u_0\rangle  &= |r\rangle , \label{eq:u_recurrence_1}\\
    |u_1\rangle  &= a_0 \tM |r\rangle  + b_0 |r\rangle , \\
    |u_{n+1}\rangle  &= a_n \tM |u_n\rangle  + b_n |u_n\rangle  - c_n |u_{n-1}\rangle . \label{eq:u_recurrence_3}
\end{align}
Thus, the $n$th Jacobi moment can be calculated as
\begin{equation}
    \mu _n = \overline{\langle r|u_n\rangle }.
\end{equation}

Practically, only a finite number of moments can be calculated. Therefore, we obtain the polynomial approximation of the form
\begin{equation}
    \rho _{\textsc{kpm}}(x) = w(x) \sum_{n=0}^{N-1} \frac{\mu _n g_n}{\hh_n} \PP_n(x),
\end{equation}
where the damping factors $g_n$ are introduced. These factors reduce the Gibbs oscillations and play a crucial role in the quality of the polynomial approximation~\cite{Silver-1994-densities-states-megadimensional}. To preserve the normalization
\begin{equation}
    \int_{-1}^1 \rho _{\textsc{kpm}}(x)\, dx = \int_{-1}^1 \rho (x)\, dx = 1,
\end{equation}
one should have $g_0 = 1$.

Practically, the averaging is performed over a finite number of random vectors $|r\rangle $. This results in slightly different Jacobi moments and the DOS, which is discussed in \ref{app:avg}.

\section{Optimal non-negative kernel}
\label{sec:optimal}

The approximated DOS can be written as
\begin{equation}
    \rho _{\textsc{kpm}}(x) = w(x) \int_{-1}^{1} K_N(x,y) \rho (y)\, dy,
\end{equation}
where $K_N(x,y)$ is the polynomial kernel of degree $N-1$ given by
\begin{equation}
	K_N(x,y) = \sum_{n=0}^{N-1} \frac{g_n}{\hh_n} \PP_n(x) \PP_n(y).   \label{eq:kernel}
\end{equation}
The DOS is essentially a non-negative quantity. Therefore, it is worthwhile to find an approximated DOS $\rho _{\textsc{kpm}}(x)$ which is also non-negative for any non-negative $\rho (x)$ \cite{Silver-1996-kernel-polynomial-approximations}. This condition is equivalent to finding the non-negative kernel
\begin{equation}
    K_N(x,y) \geq  0, \quad -1 < x, y < 1.   \label{eq:positive}
\end{equation}

There are many possibilities to find the non-negative kernel $K_N(x,y)$. One should define the objective function, which indicates the quality of the kernel. One can cosider the resolution of the kernel, which should be minimized. In this paper we consider the squared kernel resolution defined as
\begin{equation}
	Q \equiv  \frac{1}{2 \hh[\alpha +1,\beta +1]_0}\int_{-1}^1 \int_{-1}^1 (x-y)^2  K_N(x,y) w(x)w(y)\,dx\, dy,   \label{eq:Q}
\end{equation}
For Chebyshev polynomials of the first kind ($\alpha =\beta =-1/2$), the definition of $Q$ coincides with the quantity discussed in~\cite{Weisse-2006-kernel-polynomial-method}. The factor $\hh[\alpha +1,\beta +1]_0$ is chosen to make a quantitative assessment of the kernel resolution for different values of the parameters $\alpha $ and $\beta $. The resolution of the obtained kernel will be discussed in Section~\ref{sec:res} in more details.

Using the definition of the kernel given by Eq.~\eqref{eq:kernel} and the basic properties of the Jacobi polynomial, we obtain
\begin{equation}
	Q = \frac{g_0-g_1}{\alpha +\beta +2}.  \label{eq:Q_g1}
\end{equation}
Therefore, one should maximize $g_1$ under the condition \eqref{eq:positive} taking into account that $g_0=1$ due to the normalization.

\begin{figure}[t]
	\centering
	\includegraphics[scale=0.8]{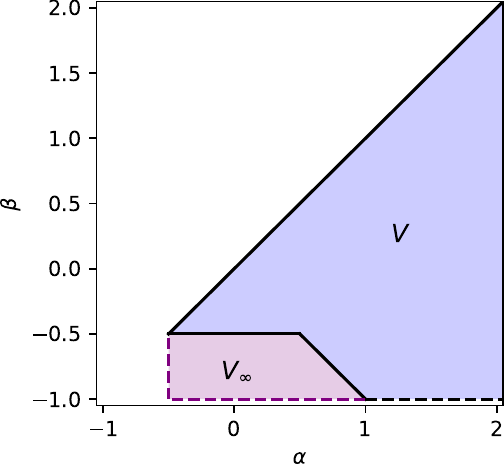}
	\caption{The region $V$ shows $(\alpha ,\beta )$ points for which $K_N(x,y)$ is non-negative for all $N$ (see text for more details). The region $V_\infty $ shows points which don't belong to $V$ but fulfill the asymptotic non-negativity condition, in which case $K_N(x,y)$ is non-negative for $N\rightarrow \infty $ (see Section~\ref{sec:asymp_opt_kernel}).}
	\label{fig:V_domain}
\end{figure}

To do so, the Gasper's theorem~[\citealp{Gasper-1972-banach-algebras-jacobi}; \citealp[Theorem 2]{Gasper-1975-positivity-special-functions}] can be used. It states that the multiplication of two Jacobi polynomials can be presented as an integral over one Jacobi polynomial,
\begin{equation}
    \PP_n(x) \PP_n(y) = \PP_n(1) \int_{-1}^1 \PP_n(z) W(x,y,z)w(z)\,dz,  \label{eq:pp}
\end{equation}
where $W(x,y,z)$ is a non-negative function for $-1 \leq  x,y,z \leq  1$ if $(\alpha ,\beta ) \in  V$, where $V$ is the set of $(\alpha ,\beta )$ such that $\alpha  \geq  \beta  > -1$, $\alpha  > -1/2$, and either $\beta  \geq  -1/2$ or $\alpha  + \beta  \geq  0$ (see Fig.~\ref{fig:V_domain}).

Below we assume that $(\alpha ,\beta ) \in  V$ and $-1 \leq  x,y,z \leq  1$ unless otherwise stated. Using Eqs.~\eqref{eq:kernel} and \eqref{eq:pp}, the kernel $K_N(x,y)$ can be presented as
\begin{equation}
    K_N(x,y) = \int_{-1}^1 K_N(z, 1) W(x, y, z)w(z)\,dz.
\end{equation}
Therefore, $K_N(x,y)\geq 0$ if $K_N(x,1) \geq  0$ because $W(x,y,z)\geq 0$. At the same time $K_N(x, 1) \geq 0$ for any non-negative $K_N(x,y)$. Therefore, the kernel $K_N(x,y)$ for $(\alpha ,\beta )\in V$ is non-negative if and only if $K_N(x,1)$ is non-negative.

Therefore, one can find the optimal non-negative polynomial $K_N(x)\equiv K_N(x,1)$. Using this polynomial, the damping factors $g_n$ can be presented as
\begin{equation} 
    g_n = \int_{-1}^{1}K_N(x)\frac{\PP_n(x)}{\PP_n(1)}w(x)\,dx \label{eq:gn_Kx}
\end{equation}
and Eq.~\eqref{eq:Q_g1} can be applied to find the resolution $Q$, which should be minimized. For simplicity, we will consider the following integrals 
\begin{equation}
    I_n = \int_{-1}^1 x^n K_N(x) w(x)\,dx. 
\end{equation}
Due to the normalization, $I_0 = g_0 = 1$. At the same time,
\begin{equation}
    g_1 = 1 - \frac{\alpha +\beta +2}{2(\alpha +1)}(1-I_1).
\end{equation}
Therefore, one should maximize $I_1$ to minimize $Q$ given by Eq.~\eqref{eq:Q_g1}. The problem to find the non-negative polynomial $K_N(x)$ of degree $N-1$ with $I_0=1$ and maximum value of $I_1$ is known as the problem of Chebyshev~[\citealp[Chapter 7.72]{Szego-1939-orthogonal-polynomials}; \citealp{Ivanov-2021-chebyshev-problem-moments}]. The optimal non-negative polynomial $K_N(x)$ is different for odd and even $N$:
\begin{gather}
    K_{2M-1}(x) = C_M^{(\alpha ,\beta )} \left(\frac{\PP_M(x)}{x - \xx_M}\right)^2,  \label{eq:Kx_odd}\\
    K_{2M}(x) = (1+x)C_M^{(\alpha ,\beta +1)}\left(\frac{\PP[\alpha ,\beta +1]_M(x)}{x - \xx[\alpha ,\beta +1]_M}\right)^2,  \label{eq:Kx_even}
\end{gather}
where $\xx_M$ is the biggest root of $\PP_M(x)$ and $C_M^{(\alpha ,\beta )}$ is the normalization constant that guarantees $I_0=1$:
\begin{equation}
    C_M^{(\alpha ,\beta )} = \frac{1-\big(\xx_M\big)^2}{(2M+\alpha +\beta +1)\hh_M}.
\end{equation}
One can show that for the optimal kernel $I_1 = \xi _N$, where $\xi _N$ is defined as
\begin{align}
	\xi _{2M-1} &= \xx_M ,\\
	\xi _{2M} &= \xx[\alpha ,\beta +1]_M,
\end{align}
for odd and even $N$ respectively.

\begin{figure}[t]
	\centering
	\includegraphics[scale=0.8]{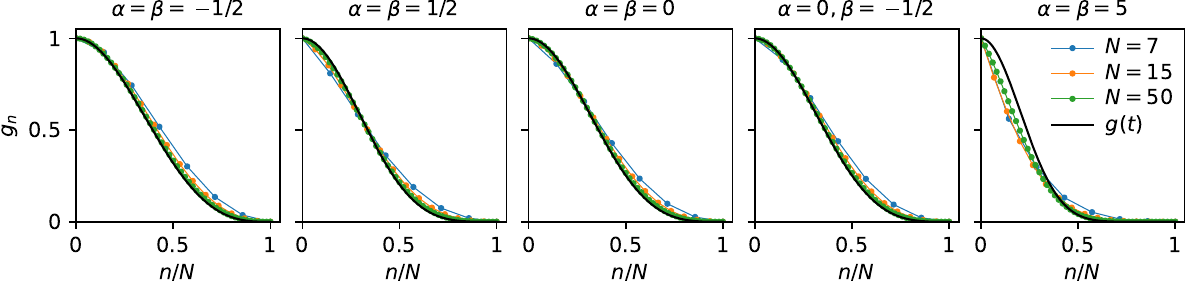}
	\caption{Damping factors $g_n$ for several small values of $N$ and their asymptotic function $g(t)$ for different values of the parameters $\alpha ,\beta $.}
	\label{fig:gs_examples}
\end{figure}

To calculate all damping factors $g_n$, one can use the integral form \eqref{eq:gn_Kx}. The first two damping factors $g_1$ and $g_2$ for the optimal kernel can be written explicitly as
\begin{gather}
	g_1 = 1 - \frac{\alpha +\beta +2}{2(\alpha +1)}(1-\xi _N), \label{eq:g_1_explicit} \\
	g_2 = 1 - (1-\xi _N) \frac{\alpha +\beta +3}{\alpha +1} \left(1 - \frac{\alpha +\beta +4}{4(\alpha +2)} \left(1 - \xi _N + \frac{1+\xi _N}{N+2+\alpha +\beta }\right)\right). \label{eq:g_2_explicit}
\end{gather}
For practical applications, one can use the Gaussian quadrature to calculate the integral in Eq.~\eqref{eq:gn_Kx} and obtain all damping factors $g_n$:
\begin{equation}\label{eq:gn_from_KN}
    g_n = \sum_{i=1}^N K_N(x_i)\frac{\PP_n(x_i)}{\PP_n(1)} w_i,
\end{equation}
where $x_i$ are zeros of $\PP_N(x)$ and $w_i$ are the corresponding Gaussian quadrature weights. There are high-performance numerical routines to obtain the values of $x_i$ and $w_i$ (for example, in the SciPy library~\cite{Virtanen-2020-scipy-fundamental-algorithms}). If one of zeros $x_i$ is close to $\xx_M$ or $\xx[\alpha ,\beta +1]_M$, there  is a possible loss of precision in Eqs.~\eqref{eq:Kx_odd}, \eqref{eq:Kx_even}, respectively. In this case, one can use the quadrature with $N+1$ points to avoid the precision loss. A minimal implementation of the KPM utilizing Jacobi polynomials with optimal damping factors computed via \eqref{eq:gn_from_KN} is available on~\cite{jacobi-kpm}.

\begin{figure}[t]
	\centering
	\includegraphics[scale=0.8]{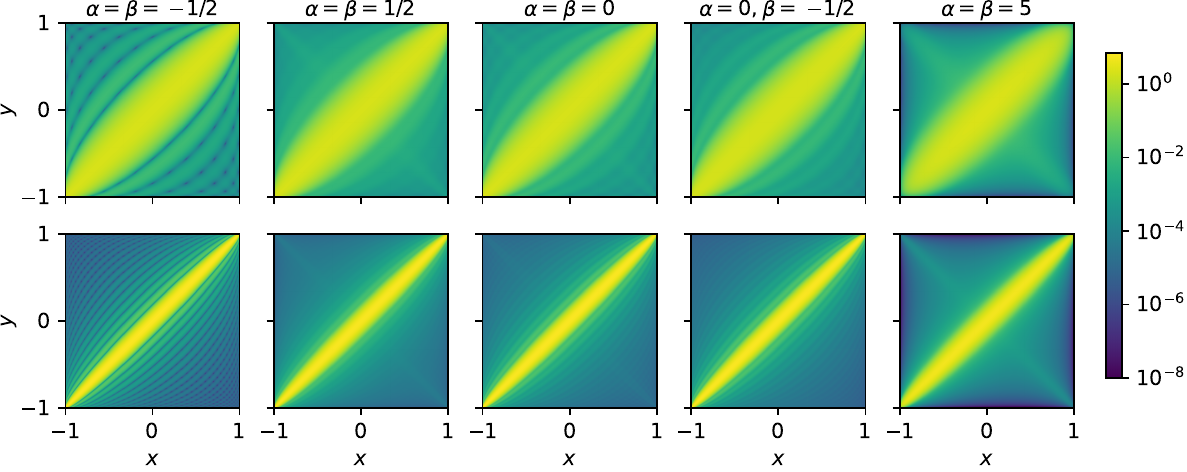}
	\caption{Heatmaps of the weighted kernel $\tilde{K}_N(x,y)$ for different values of $\alpha ,\beta $ for $N=15$ (top row) and $N=50$ (bottom row).
	}
	\label{fig:kernel_comparison}
\end{figure}

The obtained damping factors are shown in Fig.~\ref{fig:gs_examples} for different values of approximation order $N$ and the parameters $\alpha ,\beta $. One can see that the damping factors gradually decrease from $g_0=1$ to $g_N=0$. Fig.~\ref{fig:kernel_comparison} shows examples of the obtained optimal kernels for different values of $\alpha ,\beta $ and $N$. For better visual performance, a weighted kernel
\begin{equation}
	\tilde{K}_N(x,y) = K_N(x,y) \sqrt{w(x) w(y)} (1-x^2)^\frac{1}{4} (1-y^2)^\frac{1}{4}
\end{equation}
is shown. Such a multiplier makes $\tilde{K}_N(x,x)$ approximately constant regardless of $\alpha ,\beta $, as shown later in Section~\ref{sec:asymp_opt_kernel}.

\section{Resolution of the optimal kernel}
\label{sec:res}

\begin{figure}[t]
	\centering
	\includegraphics[scale=0.8]{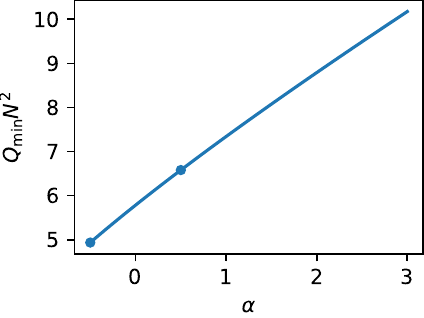}
	\caption{Asymptotic value of the scaled kernel resolution $Q_{\rm min}$ as a function of the parameter $\alpha $ for large $N$. Symbols show the values given by Eq.~\eqref{eq:sigma_exact}.}
	\label{fig:jn}
\end{figure}

For the optimal kernel, $Q$ takes its minimum value $Q_{\rm min}$. According to Eqs.~\eqref{eq:Q_g1} and \eqref{eq:g_1_explicit}, it is equal to
\begin{equation}
	Q_{\rm min} = \frac{1-\xi _N}{2(\alpha +1)}.   \label{eq:Q_opt}
\end{equation}
One can use the Mehler–Heine formula for large $N$~\cite[Chapter 7.32]{Szego-1939-orthogonal-polynomials}
\begin{equation}\label{eq:jacobi_mehler-heine}
    \PP_N\left(1 - \frac{r^2}{2N^2}\right) \simeq \left(\frac{2N}{r}\right)^\alpha  \!J_\alpha (r),
\end{equation}
where $J_\alpha (r)$ is the Bessel function. Therefore, for large $N$, the largest root $\xi _N$ both for even and odd $N$ has the asymptotic value
\begin{equation}
    \xi _N \simeq 1 - 2\frac{j_{\alpha ,1}^2}{N^2},
\end{equation}
where $j_{\alpha ,n}$ denotes the $n$th positive zero of the Bessel function $J_\alpha (x)$. Thus, the squared kernel resolution is
\begin{equation}
    Q_{\rm min} \simeq \frac{j_{\alpha ,1}^2}{1+\alpha }\frac{1}{N^2}.   \label{eq:Q_largeN}
\end{equation}
Therefore, the optimal kernel resolution $\sqrt{Q_{\rm min}}$ scales as $1/N$ for any values of the parameters $\alpha $ and $\beta $. One can see that for large $N$ the optimal squared kernel resolution $Q_{\rm min}$ depends on the parameter $\alpha $ only. This dependence is presented in Fig.~\ref{fig:jn}. The dependence on the parameter $\beta $ is implicit and given by the inequality $\beta \leq \alpha $. For $\alpha =\pm 1/2$, using $j_{-1/2,1} = \pi /2$ and $j_{1/2,1} = \pi $ we obtain
\begin{equation}
	Q_{\rm min} \simeq 
	\begin{cases}
		\dfrac{\pi ^2}{2N^2}, & \alpha  = -\dfrac{1}{2}, \\[1ex]
		\dfrac{2\pi ^2}{3N^2}, & \alpha  = \dfrac{1}{2}. \\
	\end{cases}
	\label{eq:sigma_exact}
\end{equation}
For large $\alpha $, the asymptotic behavior of the first zeros of the Bessel function is $j_{\alpha ,1} = \alpha  + \mathcal{O}(\alpha ^{1/3})$, see \cite[p.~521]{Watson-1944-treatise-theory-bessel}. Therefore, $Q_{\rm min} \simeq \alpha /N^2$ for large $\alpha $.

However, the resolution may depend on the position in the spectrum. In order to investigate it, one can consider the DOS of the form $\rho (x) = \delta (x - x_0)$. In this case the approximated DOS $\rho _{\textsc{kpm}}(x)$ is a finite-width delta-function
\begin{equation}
	\delta _{\textsc{kpm}}(x - x_0) = w(x) \int_{-1}^{1} K_N(x,y) \delta (y - x_0)\, dy = w(x) K_N(x,x_0)
\end{equation}
with the normalization $\int_{-1}^1 \delta _{\textsc{kpm}}(x - x_0)\,dx = 1$. We can find the mean and the mean square of $\delta _{\textsc{kpm}}(x - x_0)$: 
\begin{align}
	\overline{x} &\equiv  \int_{-1}^1 x \delta _{\textsc{kpm}}(x - x_0)\, dx
	= \frac{\beta -\alpha }{\alpha +\beta +2} + \left(\frac{\alpha -\beta }{\alpha +\beta +2} + x_0\right) g_1, \\
	\overline{x^2} &\equiv  \int_{-1}^1 x^2 \delta _{\textsc{kpm}}(x - x_0)\, dx
	= \frac{2+\alpha +\beta +(\alpha -\beta )^2}{(2+\alpha +\beta )(3+\alpha +\beta )} - \frac{2(\alpha -\beta )}{4+\alpha +\beta } \left(\frac{\alpha -\beta }{2+\alpha +\beta } + x_0 \right) g_1
	\notag\\
	& \hspace{4.5cm} + \left((1-x_0)^2 + \frac{4(1+\alpha )(2+\alpha )}{3+\alpha +\beta } - \frac{4(2+\alpha -x_0)(2+\alpha )}{4+\alpha +\beta }\right) g_2.
\end{align}
Using the damping factors given by Eqs. \eqref{eq:g_1_explicit}--\eqref{eq:g_2_explicit}, one can obtain the mean and the dispersion for the optimal kernel:
\begin{align}
	\overline{x} &= x_0 - Q_{\rm min}\bigl(\alpha -\beta  + (\alpha +\beta +2)x_0\bigr), \\
	\sigma ^2 &= \overline{x^2} - \overline{x}^2 \simeq 2 Q_{\rm min} \left(1-x_0^2 + \frac{4+(1+x_0)\alpha (\alpha +3) + (1-x_0)\beta (\beta +3)}{(2+\alpha )N}\right).
\end{align}
For simplicity, the dispersion is written up to the relative error of $\mathcal{O}(1/N)$. 

For large $N$, $\overline{x}$ has the expected value $\overline{x} \simeq x_0$. The dispersion is $\sigma ^2 \simeq 2 Q_{\rm min} (1-x_0^2)$ for all values of $x_0$ except a narrow interval near edge points $x_0=\pm  1$. Near edges the dispersion is $\sigma ^2 \simeq 4 Q_{\rm min} \frac{1+\alpha }{N}$ and $\sigma ^2 \simeq 4 Q_{\rm min} \frac{2+3\beta +\beta ^2}{(2+\alpha )N}$ for $x_0=1$ and $x_0=-1$, respectively. Thus, $\sigma $ scales as $1/N$ in the bulk of the spectrum and as $1/N^{3/2}$ near the edges, matching the behavior of the Jackson damping factors that are optimal for $\alpha =\beta =-1/2$~\cite{Weisse-2006-kernel-polynomial-method}.

\section{Special cases}
\label{sec:special}

For $\alpha =\pm 1/2$ and $\beta =\pm 1/2$, the Jacobi polynomials have the trigonometric form~\cite[Chapter 4.1]{Szego-1939-orthogonal-polynomials}:
\begin{align}
    \frac{1}{q_n}\PP[-\frac{1}{2},-\frac{1}{2}]_n(x) &= T_n(x) = \cos(n \varphi ),
    \label{eq:cheb1}\\
    \frac{1}{2q_{n+1}}\PP[\frac{1}{2},\frac{1}{2}]_n(x) &= U_n(x) = \dfrac{\sin((n+1)\varphi )}{\sin \varphi }, 
    \label{eq:cheb2}\\
    \frac{1}{q_n }\PP[-\frac{1}{2},\frac{1}{2}]_n(x) &= V_n(x) = \frac{\cos((n+\frac{1}{2})\varphi )}{\cos(\varphi /2)}, 
    \label{eq:cheb3}\\
    \frac{1}{q_n}\PP[\frac{1}{2},-\frac{1}{2}]_n(x) &= W_n(x) = \dfrac{\sin((n+\frac{1}{2})\varphi )}{\sin(\varphi /2)},
    \label{eq:cheb4}
\end{align}
where $\varphi =\arccos x$ and the normalization constant is
\begin{equation}
    q_n = \frac{\Gamma (n+\frac{1}{2})}{\Gamma (\frac{1}{2})\Gamma (n+1)}.
\end{equation}
The polynomials $T_n(x)$, $U_n(x)$, $V_n(x)$, and $W_n(x)$ are known as Chebyshev polynomials of the first to the fourth kind, respectively~\cite{Mason-1993-nearminimax-complex-approximation}. The trigonometric forms \eqref{eq:cheb1}--\eqref{eq:cheb2} allow finding the exact form of the damping factors $g_n(t)$ in these cases.

For $\alpha =\beta =-1/2$ and arbitrary $N$, we have:
\begin{gather}
    \xi _N = \cos\frac{\pi }{N + 1}, \\
    Q_{\rm min} = 1 - \cos\frac{\pi }{N + 1} \simeq \frac{\pi ^2}{2N^2}, \\
    g_n = \frac{N-n+1}{N+1}\cos\frac{\pi n}{N+1} +\frac{\cot\frac{\pi }{N+1}}{N+1}\sin\frac{\pi n}{N+1}.
\end{gather}
It is the known damping factors for the Jackson kernel~\cite{Weisse-2006-kernel-polynomial-method}. 

For $\alpha =1/2$, $\beta =-1/2$ and arbitrary $N$, we have
\begin{gather}
    \xi _N = \cos\frac{2\pi }{N + 2}, \\
    Q_{\rm min} = \frac{1}{3}\left(1 - \cos\frac{2\pi }{N + 2}\right) \simeq \frac{2}{3}\frac{\pi ^2}{N^2},  \\
    \begin{multlined}[b][.88\displaywidth]
        g_n = \frac{1}{2 (2n+1) (N+2)} \biggl(\cot^2\frac{\pi }{N+2}
        - \frac{1+3\cos\frac{2\pi }{N+2}}{\sin\frac{\pi }{N+2}\sin\frac{2\pi }{N+2}} \cos \frac{(2 n +1) \pi }{N+2} \\+\frac{2N-2n+3}{\sin\frac{\pi }{N+2}} \sin \frac{(2n+1) \pi }{N+2}\biggr).
    \end{multlined}
\end{gather}

For $\alpha =\beta =1/2$, the analytical expressions can be obtained only for odd $N$:
\begin{gather}
    \xi _N = \cos\frac{2\pi }{N + 3}, \\
    Q_{\rm min} = \frac{1}{3}\left(1 - \cos\frac{2\pi }{N + 3}\right) \simeq \frac{2}{3}\frac{\pi ^2}{N^2}, \\
    \begin{multlined}[b][.88\displaywidth]
        g_n = \frac{1}{2(n+1) (N+3)}\biggl(\cot^2\frac{\pi }{N+3}
        +(-1)^n \tan^2 \frac{\pi }{N+3} -  \frac{4\cos\frac{2 \pi }{N+3}}{\sin^2\frac{2 \pi }{N+3}}\cos \frac{2 (n+1) \pi }{N+3} \\
        + \frac{2(N-n+2)}{\sin\frac{2 \pi }{N+3}} \sin \frac{2 (n+1) \pi }{N+3}\biggr).
    \end{multlined}
\end{gather}    
For even $N$, the largest root $\xi _N$ as well as the other expressions could not be presented analytically. However, the value of $N$ can be decreased by one to obtain the odd number.

\section{Asymptotic behavior}

\subsection{Asymptotics of the damping factors}

In many practical applications, the value of $N$ should be large enough to obtain sufficient kernel resolution. In this case one can find damping factors $g_n$ in the limit $N\to\infty$. The polynomial $K_N(x)$ is close to 0 for all values of $x$ except a small range of $x$ near the boundary $x=1$.

For large $N$, both for even and odd $N$, using Eq.~\eqref{eq:jacobi_mehler-heine} we obtain
\begin{equation}
    K_N\left(1 - \frac{2r^2}{N^2}\right) \simeq \frac{j_{\alpha ,1}^2N^{2\alpha +2}}{2^{\alpha +\beta +1}r^{2\alpha }} \frac{J_\alpha ^2(r)}{(r^2 - j_{\alpha ,1}^2)^2}.
\end{equation}
As a result, for large $N$, the damping factors can be written as $g_n = g(n/N)$ with continuous damping function $g(t)$ defined as
\begin{equation}
    g(t) = 2j_{\alpha ,1}^2 \Gamma (\alpha +1)\int_0^\infty  \frac{J_\alpha ^2(r)}{(r^2 - j_{\alpha ,1}^2)^2} \frac{J_\alpha (2tr)}{\left(t r\right)^{\alpha }}r \,dr. \label{eq:g_cont_int}
\end{equation}

For numerical applications, the integral representation \eqref{eq:g_cont_int} is not well suited. One can show that $g(t)$ given by Eq.~\eqref{eq:g_cont_int} is equal to zero for $t\geq 1$. Therefore, one can use the Bessel-Fourier transform on the interval $0\leq t\leq 1$ and obtain the series representation of $g(t)$:
\begin{equation}
    g(t) = 2^\alpha \Gamma (\alpha +1) \sum_{n=1}^\infty  d_n^{(\alpha )} \frac{J_\alpha (t j_{\alpha ,n})}{(tj_{\alpha ,n})^\alpha },  \label{eq:g_cont}
\end{equation}
where
\begin{equation}
    d_n^{(\alpha )} = \frac{j_{\alpha ,1}^2}{J_{\alpha +1}^2(j_{\alpha ,n})}\frac{J_\alpha ^2(j_{\alpha ,n}/2)}{(j_{\alpha ,n}^2/4 - j_{\alpha ,1}^2)^2}.
\end{equation}
All coefficients $d_n^{(\alpha )}$ are non-negative and their sum is equal to one, $\sum_{n=1}^\infty  d_n^{(\alpha )} = 1$. For large $n$, the asymptotic behavior of Bessel function can be used to find the asymptotic behavior of $d_n^{(\alpha )}$:
\begin{equation}
    d_n^{(\alpha )} \simeq \frac{16j_{\alpha ,1}^2}{\pi ^4n^4} \left[1-(-1)^n\sin\left(\frac{\alpha \pi }{2} +  \frac{\pi }{4}\right)\right].
\end{equation}
The coefficients $d_n^{(\alpha )}$ decrease rapidly as $\mathcal{O}(1/n^4)$, and so not many terms in Eq.~\eqref{eq:g_cont} are required to achieve the required numerical precision for $g(t)$.

\begin{figure}[t]
	\centering
	\includegraphics[scale=0.8]{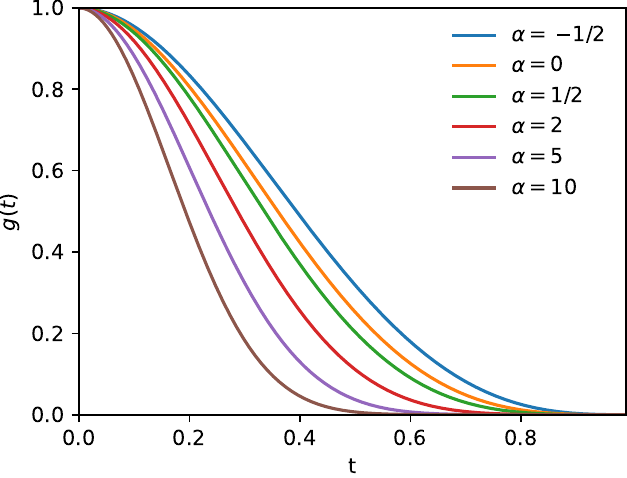}
	\caption{Asymptotic damping function $g(t)$ for different values of the parameter $\alpha $.}
	\label{fig:gt_by_alpha}
\end{figure}

The behavior of $g_n$ for different $N$ compared to $g(t)$ is presented in Fig.~\ref{fig:gs_examples}. The comparison of $g(t)$ for different values of $\alpha $ is presented in Fig.~\ref{fig:gt_by_alpha}.

\subsection{Asymptotics of the optimal kernel}
\label{sec:asymp_opt_kernel}

For large values of $N$ one can find the asymptotic structure of the kernel $K_N(x,y)$. In this case the damping factors of the form $g_n = g(n/N)$ given by Eq.~\eqref{eq:g_cont_int} can be used in Eq.~\eqref{eq:kernel}:
\begin{equation}
	K_N(x,y) = \sum_{n=0}^{N-1} \frac{g(n/N)}{\hh_n} \PP_n(x) \PP_n(y).   \label{eq:kernel_asymp1}
\end{equation}
For large $N$, one can make the asymptotic analysis of Eqs.~\eqref{eq:kernel_asymp1}. 

\begin{figure}[t]
    \centering
    \includegraphics[scale=0.8]{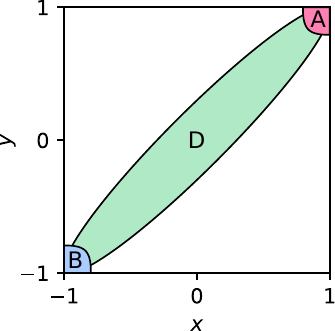}
    \caption{Regions of the kernel $K_N(x,y)$ for large $N$. The regions D, A and B are described by Eq.~\eqref{eq:kernel_asymp_diag}, \eqref{eq:kernel_asymp_A} and~\eqref{eq:kernel_asymp_B} respectively.}
    \label{fig:regions}
\end{figure}

\begin{figure}[t]
	\centering
	\includegraphics[scale=0.8]{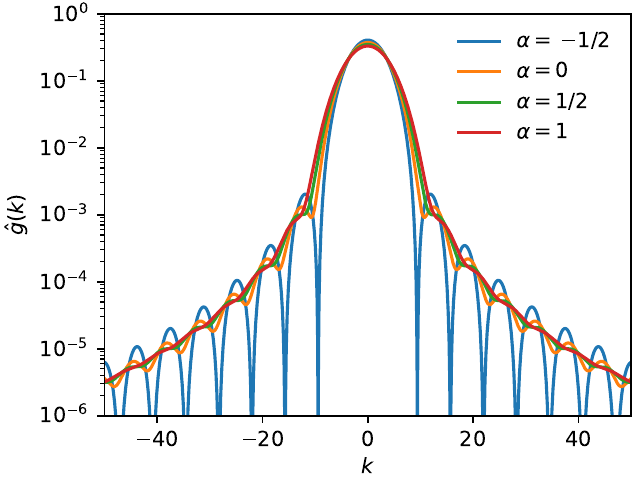}
	\caption{The cosine transform $\hat{g}(k)$ of the damping function $g(t)$ for different values of $\alpha $. It represents the kernel $K_N(x,y)$ in the diagonal region D according to Eq.~\eqref{eq:kernel_asymp_diag}.
	} 
	\label{fig:gt_cosine}
\end{figure}

For $-1 + \varepsilon  < x < 1 - \varepsilon $ and $-1 + \varepsilon  < y < 1 - \varepsilon $ with some small given $\varepsilon $, one can use Darboux's formula \cite[Chapter 7.32]{Szego-1939-orthogonal-polynomials}:
\begin{equation}
    \PP_n(x) \simeq \sqrt{\frac{2^{\alpha +\beta +1}}{n \pi  w(x)}} (1-x^2)^{-\frac{1}{4}} \cos\left[\left(n + \frac{\alpha +\beta +1}{2}\right) \arccos x - \frac{2\alpha +1}{4}\pi \right].
\end{equation}
Assuming also $\lvert\arccos x - \arccos y\rvert \lesssim 1/N$, we get the following asymptotic for the diagonal region of the kernel marked by D in Fig.~\ref{fig:regions}
\begin{equation}
	K_N(x,y) \simeq \frac{N\, \hat{g}\Bigl(N (\arccos x - \arccos y)\Bigr)}{\pi  \sqrt{w(x) w(y)}\, (1-x^2)^\frac{1}{4} (1-y^2)^\frac{1}{4}},   \label{eq:kernel_asymp_diag}
\end{equation}
where $\hat{g}(k)$ is the cosine transform of $g(t)$:
\begin{equation}\label{eq:g_cosine}
	\hat{g}(k) = \int_{0}^1 g(t) \cos(k t) dt.
\end{equation}
Thus, the kernel is given by $\hat{g}(k)$ near the diagonal $x=y$. The form of $\hat{g}(k)$ is shown in Fig.~\ref{fig:gt_cosine} for different values of the parameter $\alpha $. For large values of $k$, there is an asymptotic $\hat{g}(k) \sim k^{-4}$.

By using the definition of $g(t)$ given by Eq.~\eqref{eq:g_cont_int} and~\cite[Chapter 3.3]{Watson-1944-treatise-theory-bessel}, we obtain
\begin{equation}
	\hat{g}(k) = \frac{\Gamma (\alpha +1) \Gamma (\frac{1}{2}) j_{\alpha ,1}^2 }{2^{2\alpha -1} \Gamma (\alpha +\frac{1}{2})} \int_{|k|/2}^\infty  \frac{J_\alpha ^2(r) \,(4 r^2 - k^2)^{\alpha -\frac{1}{2}}}{(r^2 - j_{\alpha ,1}^2)^2\, r^{2\alpha -1}} dr,
\end{equation}
which is positive, and hence the asymptotic kernel in the region D is positive for all $\alpha \geq \beta \geq -1/2$.

For $x$ and $y$ close to 1 (region A in Fig.~\ref{fig:regions}), we instead apply Eq.~\eqref{eq:jacobi_mehler-heine}
to get
\begin{equation}\label{eq:region_A_asymp1}
	K_N\left(1-\frac{p^2}{2N^2}, 1-\frac{q^2}{2N^2}\right) 
	\simeq \frac{\Gamma (\alpha +1) j_{a,1}^2 N^{2\alpha +2}}{2^{\beta -\alpha -1}}
	\int_0^\infty  \frac{J_\alpha ^2(r)}{(r^2 - j_{\alpha ,1}^2)^2} \int_0^\infty  \frac{J_\alpha (p t) J_\alpha (q t)}{p^\alpha  q^\alpha }  \frac{J_\alpha (2tr)}{\left(t r\right)^{\alpha -1}} dt \, dr,
\end{equation}
where we have replaced summation by $n$ with integration by $t$, and have used the fact that $g(t)=0$ for $t>1$ to extend the integral over $t$ to infinity. Using the Sonine formula~\cite[Chapter 13.46]{Watson-1944-treatise-theory-bessel}
\begin{equation}\label{eq:bessel_triple_integral_same_arg}
	\int_0^\infty  J_\nu (a t) J_\nu (b t) J_\nu (c t) \frac{d t}{t^{\nu -1}} = \frac{2^{\nu -1} \Delta ^{2\nu -1}(a,b,c)}{(a b c)^\nu  \Gamma (\nu +\frac{1}{2}) \Gamma (\frac{1}{2})},
\end{equation}
where $\Delta (a,b,c)$ is the area of the triangle with sides $a,b,c$ (0 if there's no such triangle), for $p,q \lesssim 1$ we obtain
\begin{equation}\label{eq:kernel_asymp_A}
	K_N\left(1-\frac{p^2}{2N^2}, 1-\frac{q^2}{2N^2}\right) \simeq
	\frac{\Gamma (\alpha +1) j_{\alpha ,1}^2 N^{2\alpha +2}}{2^{\beta -\alpha } \Gamma (\alpha +\frac{1}{2}) \Gamma (\frac{1}{2})} \int_{\frac{\left|p-q\right|}{2}}^{\frac{p+q}{2}} \frac{J_\alpha ^2(r)}{(r^2 - j_{\alpha ,1}^2)^2} \frac{\Delta ^{2\alpha -1}(p, q, 2r)}{(p q r)^{2\alpha }} r dr.
\end{equation}
As $\alpha  > -1/2$, this expression is positive.

For $x,y$ close to $-1$ (region B in Fig.~\ref{fig:regions}) we obtain the following asymptotic behavior for $p,q \lesssim 1$:
\begin{equation}
	K_N\left(\frac{p^2}{2N^2}-1, \frac{q^2}{2N^2}-1\right)
	\simeq \frac{\Gamma (\alpha +1) j_{\alpha ,1}^2 N^{2\beta +2}}{2^{\alpha -\beta -1}}
	\int_{0}^\infty  \frac{J_\alpha ^2(r)}{(r^2 - j_{\alpha ,1}^2)^2} \int_0^\infty   \frac{J_\beta (p t) J_\beta (q t)}{p^\beta  q^\beta } \frac{J_\alpha (2tr)}{(t r)^{\alpha -1}} dt\,dr
\end{equation}
which becomes identical to Eq.~\eqref{eq:region_A_asymp1} when $\alpha =\beta $. For $\alpha \ne \beta $, we apply formula (1) in~\cite[Chapter 13.46]{Watson-1944-treatise-theory-bessel} to obtain 
\begin{multline}\label{eq:kernel_asymp_B} 
	K_N\left(\frac{p^2}{2N^2}-1, \frac{q^2}{2N^2}-1\right)
	\simeq \frac{2^{\beta -3\alpha +2} \Gamma (1+\alpha ) j_{\alpha ,1}^2 N^{2\beta +2}}{\Gamma \left(\frac{1}{2}\right) \Gamma \left(\beta +\frac{1}{2}\right) \Gamma (\alpha -\beta )}
	\int_0^\pi  \int_{r_0}^\infty  \frac{J_\alpha ^2(r)}{(r^2 - j_{\alpha ,1}^2)^2} \\
	\times   \left(4r^2 - p^2 - q^2 + 2 p q \cos(\phi )\right)^{\alpha -\beta -1} \sin^{2\beta }(\phi ) \frac{dr \,d\phi }{r^{2\alpha -1}},
\end{multline}
where $r_0 = \frac{1}{2} \sqrt{p^2 + q^2 - 2 p q \cos{\phi }}$. For $\alpha >\beta >-1/2$, this expression is positive.

As a result, we have shown that, in the $N\rightarrow \infty $ limit, the obtained kernel $K_N(x,y)$ is non-negative in the near-diagonal regions A, B, D. Outside these regions the kernel approaches zero for $N\rightarrow \infty $. Thereby we obtain asymptotic non-negativity of the kernel $K_N(x,y)$ for any $(\alpha ,\beta )$ such that $\alpha \geq \beta $, $\alpha >-1/2$, $\beta >-1$. It can be seen (Fig.~\ref{fig:V_domain}) that this region is larger than the region $V$ in which $K_N(x,y)$ is positive for \emph{all} $N$. On Fig.~\ref{fig:V_domain}, points outside $V$ for which asymptotic non-negativity holds are labeled $V_\infty $.

\section{Examples}
\label{sec:examples}

To demonstrate the results using the KPM based on Jacobi polynomials, we consider the DOS of the $d$-dimensional simple cubic lattice with
\begin{equation}
    M_{ij} = 
    \begin{cases}
        -1, & \text{$i$ and $j$ are neighbours}, \\
        2d, & i=j, \\
        0, & \text{otherwise}.
    \end{cases}
\end{equation}
Such matrices arise in the study of tight-binding models~\cite{Qin-2011-calculation-binding-energy}, lattice vibrations~\cite{Conyuh-2017-boson-peak-twodimensional}, and the Laplacian of graphs~\cite{Chung-1997-spectral-graph-theory, Braverman-2022-sublinear-time-spectral}, among other applications.

The $d$-dimensional simple cubic lattice has eigenvalue range $0  \leq  \varepsilon  \leq  4d$.
The edge behavior of the DOS depends on $d$, allowing to demonstrate several different cases for the choice of $\alpha ,\beta $. The asymptotics for low $\varepsilon $ is
\begin{gather}
	\rho (\varepsilon ) \sim \varepsilon ^{\frac{d}{2}-1}, \\
	\rho (4d-\varepsilon ) \sim \varepsilon ^{\frac{d}{2}-1}.
\end{gather}
We will consider the square lattice ($d=2$) and the cubic lattice ($d=3$), and compare the classical KPM based on Chebyshev polynomials of the first kind ($\alpha =\beta =-1/2$) with the polynomials that best match the edge behavior ($\alpha =\beta =d/2-1$).

\subsection{Square lattice}

\begin{figure}[t]
	\centering
	\includegraphics[scale=0.8]{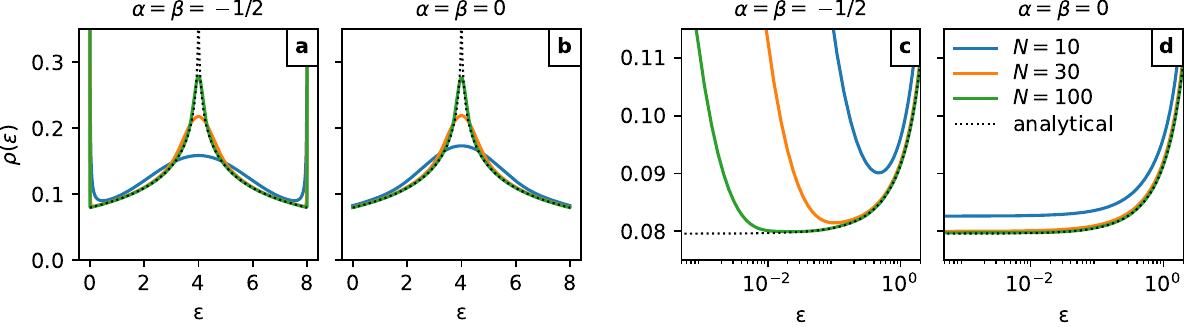}
	\caption{DOS obtained for the square lattice by the KPM for $\alpha =\beta  =- 1/2$ (a) and $\alpha =\beta =0$ (b). Panels (c) and (d) show the left tail of the DOS in detail for the corresponding values of $\alpha $ and $\beta $.
	}
	\label{fig:vdos_square}
\end{figure}

For the square lattice, the optimal parameters are $\alpha =\beta =0$, which corresponds to Legendre polynomials.
The DOS obtained by the KPM for the square lattice of size $500\times500$ is presented in Fig.~\ref{fig:vdos_square}. The figure shows a comparison of the KPM results for $\alpha =\beta =-1/2$ and for $\alpha =\beta =0$.
For stochastic averaging (see \ref{app:avg}) $R=10000$ random initial vectors were used. Figure \ref{fig:vdos_square} also shows the analytical expression of the DOS for an infinite square lattice~\cite{Hu-1981-density-states-twodimensional},
\begin{equation}
	\rho _{2 \rm D}(\varepsilon ) = \frac{1}{2\pi ^2} K\biggl[\frac{\varepsilon (8-\varepsilon )}{16}\biggr],  \label{eq:dos2d}
\end{equation}
where $K[m] = \int_0^{\pi/2} (1-m\sin^2\theta)^{-1/2} d\theta$ is the complete elliptic integral of the first kind.

The approximation accuracy of the DOS obtained by the KPM in the center of the spectrum increases with degree $N$ for both families of polynomials. However, the approximation near the edges diverges for $\alpha =\beta =-1/2$, and the quality of such an approximation increases only slowly with increasing $N$. Meanwhile, for $\alpha =\beta =0$, the correct edge behavior (a constant one) is achieved for arbitrary $N$, and the convergence is much faster.

\subsection{Cubic lattice}

\begin{figure}[t]
	\centering
	\includegraphics[scale=0.8]{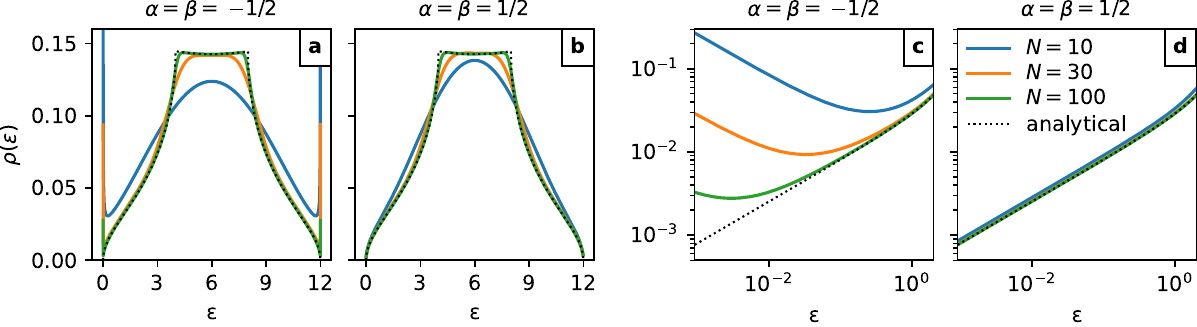}
	\caption{DOS obtained for the cubic lattice by the KPM for $\alpha =\beta =-1/2$ (a) and $\alpha =\beta =1/2$ (b). Panels (c) and (d) show the left tail of the DOS in detail for the corresponding values of $\alpha $ and $\beta $.
	}
	\label{fig:vdos_cubic}
\end{figure}

For the cubic lattice, the optimal parameters are $\alpha =\beta =1/2$, which corresponds to the Chebyshev polynomials of the second kind $U_n(x)$. The DOS obtained by the KPM for the cubic lattice of size $75\times 75\times 75$, averaged over $R=2000$ random initial vectors, is presented in Fig.~\ref{fig:vdos_cubic}. As with the previous example, the $T_n(x)$ polynomials cause the approximated DOS to go to infinity at the edges, whereas the $U_n(x)$ polynomials always achieve the correct edge behavior. Figure~\ref{fig:vdos_cubic} also shows the analytical expression of DOS for an infinite cubic lattice~\cite{Joyce-1998-cubic-modular-transformation, Conyuh-2021-random-matrix-approach}
\begin{equation}
    \rho _{3 \rm D}(\varepsilon ) = -\frac{1}{2\pi }\Im W_s\left(\frac{\varepsilon }{2} - 3 + 0i\right),  \label{eq:dos3d}
\end{equation}
where $W_s(z)$ is the third Watson integral
\begin{gather}
    W_s(z) = \frac{4(1-9\xi ^4)}{\pi^2z(1-\xi )^3(1+3\xi )} K^2\biggl[\frac{16 \xi ^3}{(1-\xi )^3(1+3\xi )}\biggr],  \label{eq:Watson} 
    \\
    \xi  = \sqrt{\frac{\sqrt{z}-\sqrt{z-1/z}}{\sqrt{z}+\sqrt{z-9/z}}}.  \label{eq:Watson2}
\end{gather}
Equation \eqref{eq:Watson2} is written in such a way that the primary branch of each square root is used. To calculate Eq.~\eqref{eq:dos3d}, the term $+0i$ can be omitted as long as the branch $\sqrt{-\varepsilon } = i\sqrt{\varepsilon }$ is used for real $\varepsilon >0$.

\section{Conclusion}

The kernel polynomial method was generalized to arbitrary Jacobi polynomials $\PP_n(x)$, allowing polynomial approximations to be calculated that have an arbitrary power law behavior at the edges of the spectrum. The optimal damping factors $g_n$ were obtained, which minimize the squared spectrum resolution $Q$ and ensure non-negativity of the polynomial approximation.

For some special cases $\alpha =\pm 1/2, \beta =\pm 1/2$, optimal $g_n$ were obtained in explicit trigonometric form presented in Section~\ref{sec:special}. In the case $\alpha =\beta =-1/2$, damping factors $g_n$ coincided with the known Jackson damping factors. For other cases, the optimal damping factors $g_n$ can be calculated with arbitrary precision by use of Gauss-Jacobi quadrature, Eq.~\eqref{eq:gn_from_KN}. For large $N$, the asymptotic function $g(t)$ given by \eqref{eq:g_cont} can be used.

The obtained damping factors ensure non-negativity of the polynomial approximation for the parameters $(\alpha ,\beta )\in V$ (see Fig.~\ref{fig:V_domain}) for arbitrary approximation order $N$. The case $\alpha <\beta $ can be mapped to the case $\alpha >\beta $ by swapping $\varepsilon _\text{max}$ and $\varepsilon _\text{min}$ in Eq.~\eqref{eq:M_tilde} along with swapping $\alpha $ and $\beta $. For $(\alpha ,\beta )\in V_\infty $ the non-negativity condition is satisfied asymptotically for $N\to \infty $. The parameters $\alpha ,\beta $ have a natural bound $\alpha ,\beta >-1$ since the DOS should be normalized. Thus, the non-negativity condition is satisfied for all possible values of $\alpha $ and $\beta $ (at least, in asymptotic sense) except a small region $\alpha  < \beta  \leq  -1/2$ or $\beta  < \alpha  \leq  -1/2$, which rarely appears in practical applications. 

The main properties of $g_n$ and their asymptotic function $g(t)$ were obtained. Using them, the asymptotic expressions for $K_N(x,y)$ for large $N$ were obtained. The resolution of the kernel scales as $1/N$ with a prefactor which depends on $\alpha $. Near edges the resolution scales as $N^{-3/2}$. The obtained damping factors $g_n$ can be used as a generalization of the Jackson damping factors for arbitrary Jacobi polynomials.

As an example, the implemented method was used to compute the DOS of cubic lattices in 2 and 3 dimensions, for which an analytic solution exists to compare with. The results show a fast convergence to theoretical results including the spectrum edges.

Existing KPM implementations (e.g. \cite{Groth-2014-kwant-software-packagea, Moldovan-2020-pybinding-v0-python}) can also be straightforwardly modified to support Jacobi polynomials, by providing appropriate recurrence relations \eqref{eq:u_recurrence_1}--\eqref{eq:u_recurrence_3} and calculating the optimal $g_n$ via Eq.~\eqref{eq:gn_from_KN}. We also provide a reference minimal implementation of the KPM based on Jacobi polynomials \cite{jacobi-kpm}.

\section{Acknowledgments}

The financial support of the Russian Science Foundation under the grant No. 22-72-10083 is gratefully acknowledged.

\appendix

\section{Stochastic averaging}
\label{app:avg}

Using the eigenvectors $|v_j\rangle $, the Jacobi moment given by Eq.~\eqref{eq:Jmom_rnd} can be presented as
\begin{equation}
    \mu _k = \sum_{i,j=1}^\mathcal{N} \overline{\langle r|v_i\rangle \langle v_j|r\rangle }\left\langle v_i\middle| \PP_n(\tM)\middle|v_j\right\rangle ,
\end{equation}
where the overline denotes averaging over the random $\mathcal{N}$-dimensional unit vector $|r\rangle $.
Due to the property
\begin{equation}
    \overline{\langle r|v_i\rangle \langle v_j|r\rangle } = \frac{1}{\mathcal{N}}\delta _{ij},
\end{equation} 
the Jacobi moments defined by Eq.~\eqref{eq:Jmom_rnd} and Eq.~\eqref{eq:Jmom_exact} are the same.

Practically, the averaging in Eq.~\eqref{eq:Jmom_rnd} is performed over a finite number of random unit vectors $|r_1\rangle , \ldots, |r_R\rangle $. In this case, the approximation of $n$th Jacobi moment is
\begin{equation}
    m_n = \frac{1}{R} \sum_{k=0}^{R-1}\left\langle r_k\middle|\PP_n(\tM)\middle|r_k\right\rangle .
\end{equation}
One can see that $\mu _n = \lim_{R\to\infty} m_n$. For finite $R$, the DOS calculated using KPM is
\begin{equation}
    \rho _{\textsc{kpm}}^{(R)}(x) = w(x) \sum_{n=0}^{N-1} \frac{m_n g_n}{\hh_n} \PP_n(x).
\end{equation}
It is the polynomial approximation of the weighted DOS
\begin{equation}
    \rho ^{(R)}(x) = \sum_{j=1}^\mathcal{N} \varkappa _j \delta (x - x_j),
\end{equation}
where
\begin{equation}
    \varkappa _j = \frac{1}{R}\sum_{k=1}^R \Big|\langle r_k|v_j\rangle \Big|^2.
\end{equation}
One can see that $\sum_j \varkappa _j = 1$ for arbitrary $R$. Therefore,
\begin{equation}
    \int_{-1}^1 \rho ^{(R)}(x)\, dx = 1,
\end{equation}
which means the same normalization for the approximated DOS:
\begin{equation}
    \int_{-1}^1 \rho _{\textsc{kpm}}^{(R)}(x)\, dx = 1
\end{equation}
In the limit $R\to\infty$, all weights $\varkappa _j$ are equal to $1/\mathcal{N}$, which means
\begin{align}
    \rho (x) &= \lim_{R\to\infty} \rho ^{(R)}(x), \\
    \rho _{\textsc{kpm}}(x) &= \lim_{R\to\infty} \rho _{\textsc{kpm}}^{(R)}(x).
\end{align}

Since $\rho _{\textsc{kpm}}^{(R)}(x)$ is the polynomial approximation of $\rho ^{(R)}(x)$, one can apply the properties of the polynomial kernel $K_N(x,y)$ for finite number $R$. Therefore, for damping factors $g_n$ defined in Section~\ref{sec:optimal}, $\rho _{\textsc{kpm}}^{(R)}(x)$ is non-negative for $-1 < x < 1$ for arbitrary numbers $N$ and $R$.

\bibliographystyle{elsarticle-num} 
\bibliography{KPM.bib}

\end{document}